\documentclass[11pt,leqno]{article} 
\usepackage{graphics}
\newtheorem{thm}{Theorem}[section]
\newtheorem{lma}{Lemma}[section]

\newcommand{\beqa}{\begin{eqnarray}}
\newcommand{\eeqa}{\end{eqnarray}}

\newcommand{\pf}{\noindent {\bf Proof:} $\s$ }
\newcommand{\epf}{ \hfill$\diamondsuit$ \medskip}

\newcommand{\beq}{\begin{equation}}
\newcommand{\eeq}{\end{equation}}
\newcommand{\lbl}{\label}
\newcommand{\s}{\; \;}

\newcommand{\ep}{\epsilon}

\newcommand{\ra}{\rightarrow}
\newcommand{\al}{\alpha}

\title{A periodic model for the dynamics of cell volume}

\author{
Philip Korman   \\ 
Department of Mathematical Sciences \\ 
University of Cincinnati \\ 
Cincinnati Ohio 45221-0025 \\
kormanp@ucmail.uc.edu
}

\date{}

\begin{document}

\maketitle
\begin{abstract} 
We prove the existence and uniqueness of positive periodic solution for a model describing the dynamics of cell volume flux, introduced in Julio A. Hernandez  \cite{H}. We also show that the periodic solution is a global attractor. Our results confirm the conjectures made in an interesting recent book of P.J. Torres \cite{T}.
 \end{abstract}

\begin{flushleft}
\end{flushleft}

\begin{flushleft}
\end{flushleft}

\section{Introduction}
\setcounter{equation}{0}
\setcounter{thm}{0}
\setcounter{lma}{0}

In \cite{H} Julio A. Hernandez proposed a general model for describing the dynamics of cell volume related to transport of water and solute across the cell membrane. The interdependence between the mass of solute $x(t)$ and water volume $y(t)$ is governed by the system 
\beqa
\lbl{0}
& x'=\alpha(t)-\beta \frac{x}{y} \\
& y'=-\gamma(t)+\sigma  \frac{x}{y}+\frac{\epsilon}{y}\,. \nonumber
\eeqa
Here $\alpha(t)>0$ represents the sources of solute, $\gamma(t)>0$ is related to decrease of the water volume, while the positive constants $\beta$, $\sigma$ and $\epsilon$ are biological interaction coefficients. As was shown in \cite{H}, this model unifies a number of other solute-solvent flux models, previously considered in the biological literature. This model is also described in detail in a recent book of P.J. Torres \cite{T}, see also P.J. Torres \cite{T1},  and J.D. Benson et al \cite{c}. As explained in \cite{T}, it is natural to assume that $\alpha(t)$ and $\gamma(t)$ are periodic functions, which is  related to circadian clocks. All of the coefficients in (\ref{0}) are assumed to be positive, and we are looking for positive and periodic solution, with components $x(t)$ and $y(t)$. It is not hard to state a necessary condition for the existence of  periodic solution (see  (\ref{5}) below). P.J. Torres \cite{T} proved that the necessary condition is also sufficient. Moreover, he conjectured that the periodic solution is unique and asymptotically stable.
\medskip

The system (\ref{0}) is of {\em cooperating type}, or a {\em monotone system}, see M.W. Hirsh \cite{H1}, H.L. Smith \cite{S}. We show that this fact allows one to apply the method of monotone iterations, where the trick is in constructing the appropriate supersolutions. We thus obtain an alternative proof of existence of solutions.  Moreover, the method of monotone iterations shows the existence of maximal and minimal solutions, from which we conclude the uniqueness. We also show that the periodic solution attracts all other positive solutions, proving the conjectures of P.J. Torres \cite{T}. The more general model in P.J. Torres \cite{T1} is also of cooperating type.

\section{The results}
\setcounter{equation}{0}
\setcounter{thm}{0}
\setcounter{lma}{0}
Any function $b(t) \in C[0,p]$, may be  decomposed as $b(t)=\bar b+ \tilde b(t)$, with $\bar b=\frac{1}{p} \int _0^p  b(s) \, ds$, and $\int _0^p  \tilde b(s) \, ds=0$.
The following lemma is proved by a direct integration.

\begin{lma}\lbl{lma:0}
Consider the equation 
\[
y'=b(t) \,,
\]
with  a  continuous $p$-periodic function $b(t)$. This equation has a $p$-periodic  solution if and only if $\int_0^p b(t) \, dt=0$.

\end{lma}
\begin{lma}\lbl{lma:1}
Consider the equation
\beq
\lbl{20}
y'+ay=b(t) \,,
\eeq
with a positive constant $a$, and a  continuous $p$-periodic function $b(t)$. The problem (\ref{20}) has a  solution of period $p$. This solution is unique, and it attracts all other solutions of (\ref{20}), as $t \ra \infty$. If $b(t)$ is positive, so is the $p$-periodic  solution.
\end{lma}

\pf
The general solution is
\[
y(t)=y_0e^{-at}+e^{-at} \int _0^t e^{as} b(s) \, ds \,.
\]
This solution is $p$-periodic, provided that $y(p)=y(0)=y_0$, which gives
\beq
\lbl{21}
y_0=\frac{1}{e^{ap}-1} \int _0^p e^{as} b(s) \, ds \,.
\eeq
If $z(t)$ is another solution of (\ref{20}), their  difference $w(t)=z(t)-y(t)$ is equal to  $w(t)=e^{-at} w(0) \ra 0$, as $t \ra \infty$, 
 proving that all solutions tend to the periodic solution $y(t)$. In particular, this implies that the periodic solution is unique.
\epf

Observe that in case $\bar b \ne 0$, the $p$-periodic solution tends to infinity, as $a \ra 0$.
\begin{lma}\lbl{lma:20}
Let $y(t)$ be the $p$-periodic solution of (\ref{20}). Then
\[
\lim _{a \ra 0} a \, y(t)=\bar b \,.
\]
\end{lma}

\pf
We have
\[
\lim _{a \ra 0} a \,y(t)=\lim _{a \ra 0} \frac{a}{e^{ap}-1} \int _0^p e^{as} b(s) \, ds=\frac{1}{p} \int _0^p  b(s) \, ds=\bar b \,,
\]
uniformly in $t$.
\epf

Consider a system
\beqa
\lbl{1}
& x'(t)=f(t,x(t),y(t)) \\
& y'(t)=g(t,x(t),y(t)) \,, \nonumber
\eeqa
where the given differentiable functions $f(t,x,y)$ and $g(t,x,y)$ are assumed to be $p$-periodic in $t$, for all $(x,y)$. We say that a pair of  $p$-periodic  differentiable functions $(a(t),b(t))$ forms a {\em subsolution pair} if for all $t$
\beqa
\lbl{2}
& a'(t) \leq f(t,a(t),b(t)) \\
& b'(t)\leq g(t,a(t),b(t)) \,. \nonumber
\eeqa
A {\em supersolution pair} $(A(t),B(t))$ is defined by reversing the inequalities in (\ref{2}). We say that sub- and supersolution pairs are {\em ordered} if $a(t)<A(t)$ and $b(t)<B(t)$ for all $t$. 

\begin{thm}\lbl{thm:1}
Assume that the problem (\ref{1}) has   ordered  sub- and supersolution pairs $(a(t),b(t))$ and   $(A(t),B(t))$. Assume that the system  (\ref{1}) is of cooperating type, i.e., for all $t \in R$, $x \in (a(t),A(t))$, $y \in (b(t),B(t))$ we have 
\beq
\lbl{3}
f_y(t,x,y) \geq 0, \s \mbox{and} \s g_x(t,x,y) \geq 0  \,.
\eeq
Then the problem (\ref{1}) has a $p$-periodic  solution, satisfying $a(t)<x(t)<A(t)$,  $b(t)<x(t)<B(t)$, for all $t$. Moreover, one can construct two monotone sequences of  $p$-periodic approximations $(x_n(t),y_n(t))$ and $(X_n(t),Y_n(t))$, which converge respectively to the minimal solution $(\underline x, \underline y)$ and to the maximal solution $(\overline x(t), \overline y(t))$. Furthermore, any solution of (\ref{1}), with the initial data satisfying $a(0)<x(0)<A(0)$ and  $b(0)<y(0)<B(0)$, converges to the product of the strips $(\underline x(t), \overline x(t)) \times (\underline y(t), \overline y(t))$.
\end{thm}

\pf
Beginning with $(x_0,y_0)=(a(t),b(t))$, we construct the sequence $(x_n(t),y_n(t))$ by calculating the $p$-periodic  solutions of the following equations 
\beqa  
\lbl{*}
& x_n'+Mx_n=Mx_{n-1}+f(t,x_{n-1},y_{n-1}) \,, \s n=1,2,\ldots \\  
& y_n'+My_n=My_{n-1}+g(t,x_{n-1},y_{n-1}) \,, \s n=1,2,\ldots \,. \nonumber
\eeqa
Here the constant $M>0$ is chosen so that the functions $Mx+f(t,x,y)$ and $My+g(t,x,y)$ are both increasing in $x$ and $y$, for $x \in [a(t),A(t)]$, $y \in [b(t),B(t)]$, for all $t \in [0,p]$. Since these intervals are compact, such $M$ exists. The sequence $(X_n(t),Y_n(t))$ is constructed similarly, beginning with $(X_0,Y_0)=(A(t),B(t))$. Using Lemma \ref{lma:1}, a standard argument shows that (componentwise)
\[
(a(t),b(t))<(x_1,y_1)<\cdots <(x_n,y_n)<\cdots <(X_n,Y_n)<\cdots 
\]
\[
<(X_1,Y_1)<(A(t),B(t)) \,.
\]

\noindent
It follows that both sequences $\{x_n(t) \}$ and $\{y_n(t) \}$ converge. Define $\underline x(t) =\lim _{n \ra \infty} x_n(t)$, and $\underline y(t) =\lim _{n \ra \infty} y_n(t)$. Passing to the limit in the integral version of (\ref{*}), we see that $(\underline x(t), \underline y(t))$ is a $p$-periodic solution of (\ref{1}). By a standard argument, this is the minimal solution, i.e., any other $p$-periodic solution of (\ref{1}) satisfies $\underline x(t) \leq x(t)$, and $\underline y(t) \leq y(t)$, for all $t$.

\medskip

To show that all solutions are attracted to the interval between the minimal and the maximal solutions, we proceed similarly to E.N. Dancer \cite{d}. We shall show that for any solution $(x(t),y(t))$ of the system (\ref{1}) with the initial data satisfying 
\[
a(0)<x(0)<A(0) \,, \s b(0)<y(0)<B(0) \,,
\]
we have
\[
x_n(t)<x(t)<X_n(t) \,, \s y_n(t)<y(t)<Y_n(t) \,,
\]
for any $n$, provided that $t$ is sufficiently large. We prove next that $x_n(t)<x(t)$ and $y_n(t)<y(t)$, with the other inequalities being similar. 
\medskip

We claim that $x(t)>a(t)$ and $y(t)>b(t)$, for all $t$. Letting $w(t)=x(t)-a(t)$ and $z(t)=y(t)-b(t)$, we see from (\ref{1}) and (\ref{2}) that
\beqa \nonumber
& w' \geq p(t)w+q(t)z \,, \s w(0)>0 \\ \nonumber
&  z' \geq r(t)w+s(t)z \,, \s z(0)>0 \,. \nonumber
\eeqa
where $p(t)=\int _0^1 f_x(t,\theta x(t)+(1-\theta)a(t), \theta y(t)+(1-\theta)b(t)) \, d \theta $, with similar expressions for $q(t) \geq 0$, $r(t) \geq 0$, and $s(t) $. If we define $\mu (t)=e^{-\int p(t) \, dt}$ and $\nu (t)=e^{-\int s(t) \, dt}$, these inequalities imply that $\left( \mu(t)w \right)' \geq \mu(t)q(t) z$, and $\left( \nu(t)z \right)' \geq \nu(t)r(t) w$. We see that $\mu (t)w(t)$ and $\nu(t) z(t)$ are positive and increasing functions, and the claim follows.
\medskip

Define the functions $\xi (t)$ and $\eta (t)$ by integrating  the following equations
\beqa
\lbl{30}
& \xi '+M\xi=Ma(t)+f(t,a(t),b(t)) \,, \s \xi (0)=a(0) \\
& \eta '+M \eta=Mb(t)+g(t,a(t),b(t)) \,, \s \eta (0)=b(0) \,, \nonumber
\eeqa
and rewrite (\ref{1}) as 
\beqa
\lbl{30.1}
& x'(t)+Mx(t)=Mx(t)+f(t,x(t),y(t)) \\
& y'(t)+My(t)=My(t)+g(t,x(t),y(t)) \,. \nonumber
\eeqa
By the above claim, the right hand sides in (\ref{30.1}) are pointwise greater than the ones in (\ref{30}). It follows that for all $t>0$, $x(t)>\xi (t)$ and $y(t)>\eta (t)$, and moreover, $x(t)-\xi (t)>x(0)-a(0)$ and $y(t)-\eta (t)>y(0)-b(0)$. By the definition of $(x_1,y_1)$ and Lemma \ref{lma:1}, $\xi (t) \ra x_1(t)$ and $\eta (t) \ra y_1(t)$, as $t \ra \infty$. Hence, at some $t_1>0$, $x(t_1)>x_1(t_1)$ and $y(t_1)>y_1(t_1)$. We now take $t_1$ as a new origin, and repeat this argument, showing that $x(t)>x_1(t)$ and $y(t)>y_1(t)$ for all $t>t_1$, and at some $t_2>t_1$, $x(t_2)>x_2(t_2)$ and $y(t_2)>y_2(t_2)$, and so on. Observe that $(x_n,y_n)$ is a subsolution pair, for each $n$.
\epf

We now consider the system 
\beqa
\lbl{4}
& x'=\alpha(t)-\beta \frac{x}{y} \\
& y'=-\gamma(t)+\sigma  \frac{x}{y}+\frac{\epsilon}{y}\,. \nonumber
\eeqa
Here $\al (t)$ and $\gamma(t)$ are positive $p$-periodic  functions, $\beta$, $\sigma$ and $\epsilon$ are positive constants.

\begin{thm}\lbl{thm:2}
The condition
\beq
\lbl{5}
\beta \bar \gamma -\sigma \bar \alpha>0
\eeq
is necessary and sufficient for the existence of positive $p$-periodic  solution of (\ref{4}). In case (\ref{5}) holds, the positive solution is unique, and it attracts all other positive  solutions of (\ref{4}), as $t \ra \infty$.
\end{thm}

\pf
Multiplying the first equation in (\ref{4}) by $\sigma$, the second one by $\beta$, adding the results and integrating, we get
\beq
\lbl{5.1}
\beta \bar \gamma -\sigma \bar \alpha=\frac{\epsilon \beta}{p} \int_0^p \frac{1}{y(t)} \, dt>0 \,,
\eeq
proving the necessity of the condition (\ref{5}).
\medskip

We now apply the Theorem \ref{thm:1}. For a subsolution pair, we consider $(p,q)$, where $p$ and $q$ are two small constants, with $\beta \frac{p}{q}<\min _{t \in R} \al (t)$. The conditions (\ref{2}) require
\beqa \nonumber
& 0<\alpha(t)-\beta \frac{p}{q} \\ \nonumber
& 0<-\gamma(t)+\sigma  \frac{p}{q}+\frac{\epsilon}{q}\,, \nonumber
\eeqa
which holds for $p$ and $q$ sufficiently small.
\medskip

We now construct a supersolution pair $(A(t),B(t))$. We choose $A(t)$ to be the positive $p$-periodic  solution of 
\[
A'(t)=\alpha(t)+\theta -\beta \frac{A(t)}{M} \,,
\]
with the constants $\theta >0$ small, and $M>0$ large to be fixed later.
Let $y_0(t)$ be the $p$-periodic  solution of
\[
y_0'(t)=-\tilde \gamma (t) \,.
\]
We set $B(t)=M+y_0(t)$. The conditions (\ref{2}) for supersolutions become
\beqa \nonumber
& A'=\alpha(t)+\theta -\beta \frac{A}{M}>\alpha(t) -\beta \frac{A}{M+y_0} \\ \nonumber
& B'=- \tilde \gamma(t)>- \tilde \gamma(t)- \bar \gamma +\sigma  \frac{A}{M+y_0}+\frac{\epsilon}{M+y_0}\,. \nonumber
\eeqa
The first of these inequalities simplifies to read
\beq
\lbl{6}
\theta > \beta \frac{A}{M} -\beta \frac{A}{M+y_0}=\beta \frac{Ay_0}{M(M+y_0)} \,.
\eeq
The second inequality requires
\beq
\lbl{7}
0>- \bar \gamma+\sigma \frac{A}{M}+\sigma  \left(\frac{A}{M+y_0}-\frac{A}{M} \right)+\frac{\epsilon}{M+y_0} \,.
\eeq
By Lemma \ref{lma:20}, $\frac{A}{M} \ra \frac{1}{\beta} \bar \al + \frac{1}{\beta} \theta$, as $M \ra \infty$. We now fix $\theta>0$ so small that
\[
0>- \bar \gamma+\sigma \frac{1}{\beta} \bar \al +\frac{\sigma}{\beta} \theta \,,
\]
and then choosing $M$ sufficiently large, we can satisfy both the inequalities (\ref{6}) and (\ref{7}), and the Theorem \ref{thm:1} applies, proving the existence. 
\medskip

By (\ref{5.1}), the $y$ component of any positive $p$-periodic  solution is equal to that of the maximal solution. But then from the first equation in (\ref{4}), the first components also coincide, proving the uniqueness.
\medskip

Observe that with our construction, we can make the supersolutions arbitrarily large, and the subsolutions arbitrarily  small. By the Theorem \ref{thm:1}, the unique $p$-periodic solution is then a global attractor for all positive solutions.
\epf

We conclude with a numerical example. We used {\em Mathematica} to solve   the problem (\ref{4}) with $\al (t)=2+\sin (2 \pi t)$, $\gamma (t)=1+ \cos ^2  (2 \pi t)$, $\beta=2$, $\sigma=1$, $\ep =  0.2$, and the initial conditions $x(0)=1$, $y(0)=0.4$. The Theorem \ref{thm:1} applies. In Figure $1$ one can see quick convergence to the unique periodic solution of period $1$. We saw similar results for all other initial conditions, and all other systems, that we tried.

\begin{figure}
\begin{center}
\scalebox{.84}{\includegraphics{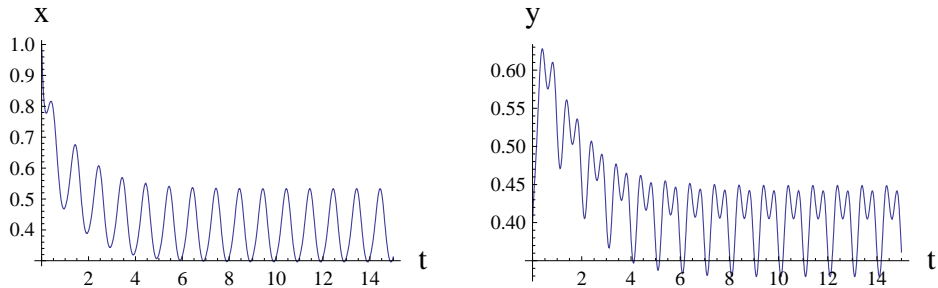}}
\end{center}
\caption{ A positive  solution of   the problem (\ref{4}) approaching the periodic solution (of period $1$) }
\end{figure}
\medskip

\noindent
{\bf Acknowledgment:} I wish to thank the referee for useful comments.

\end{document}